\documentclass[11pt]{article}
\usepackage{amscd}
\usepackage{amsmath}
\usepackage{latexsym}
\usepackage{amsfonts}
\usepackage{amssymb}
\usepackage{amsthm}
\usepackage{graphicx}

 \newtheorem{proposition}{Proposition}[section]
 \newtheorem{definition}{Definition}[section]
 \newtheorem{lemma}{Lemma}[section]
 \newtheorem{theorem}{Theorem}[section]
 \newtheorem{corollary}{Corollary}[section]
 \newtheorem{remark}{Remark}[section]

\makeatletter
\newcommand{\Extend}[5]{\ext@arrow0099{\arrowfill@#1#2#3}{#4}{#5}}
\makeatother

\begin{document}

 \title{ Global well-posedness and scattering for the energy-critical, defocusing Hartree equation for radial data}
 \author{{Changxing Miao,\ \ Guixiang Xu,\ \ and \ Lifeng Zhao }\\
         {\small Institute of Applied Physics and Computational Mathematics}\\
         {\small P. O. Box 8009,\ Beijing,\ China,\ 100088}\\
         {\small (miao\_changxing@iapcm.ac.cn, \ xu\_guixiang@iapcm.ac.cn, zhao\_lifeng@iapcm.ac.cn ) }\\
         \date{}
        }
\maketitle

\begin{abstract} We consider the defocusing, $\dot{H}^1$-critical Hartree
equation for the radial data in all dimensions $(n\geq 5)$. We
show the global well-posedness and scattering results in the
energy space. The new ingredient in this paper is that we first
take advantage of the term $\displaystyle - \int_{I}\int_{|x|\leq
A|I|^{1/2}}|u|^{2}\Delta \Big(\frac{1}{|x|}\Big)dxdt$ in the
localized Morawetz identity to rule out the possibility of energy
concentration, instead of the classical Morawetz estimate
dependent of the nonlinearity.
\end{abstract}

 \begin{center}
 \begin{minipage}{120mm}
 \end{minipage}
 \end{center}

\section{Introduction}
\setcounter{section}{1}\setcounter{equation}{0} In this paper, we
study the Cauchy problem for the Hartree equation
\begin{equation} \label{equ1}
\left\{ \aligned
    iu_t +  \Delta u  & = f(u), \quad  \text{in}\  \mathbb{R}^n \times \mathbb{R}, \quad n\geq 5,\\
     u(0)&=\varphi(x), \quad \text{in} \ \mathbb{R}^n.
\endaligned
\right.
\end{equation}
Here $f(u)=\big(V* |u|^2 \big)u$ is a nonlinear function of Hartree
type for $V(x)=|x|^{-\gamma}, 0<\gamma<n$, where $*$ denotes the
convolution in $\mathbb{R}^{n}$. In practice, we use the integral
formula of $(\ref{equ1})$
\begin{equation}\label{intequa1}
u(t)=U(t)\varphi -i \int^{t}_{0} U(t-s)f(u(s))ds,
\end{equation}
where $U(t)=e^{it\Delta}$.

If the solution $u$ of $(\ref{equ1})$ has sufficient smoothness and
decay at infinity, it satisfies two conservation laws :
\begin{equation}\label{cl}
\aligned
& M(u(t)) =\big\|u(t)\big\|_{L^2} = \big\|\varphi\big\|_{L^2},\\
E(u(t))= \frac12 & \big\|\nabla u(t)\big\|^2_{L^2}+\frac{1}{4} \int
\int \frac{1}{|x-y|^{\gamma}} |u(t,x)|^2 |u(t,y)|^2\ dxdy
=E(\varphi).
\endaligned
\end{equation}
As explained in \cite{GiV00}, the energy is also conserved for the
energy solutions $u\in C^{0}_t(\mathbb{R}, H^1)$.

From the viewpoint of the fractional integral, we rewrite the
equation $(\ref{equ1})$ as
\begin{equation*}
    iu_t +  \Delta u   =\big( (-\Delta)^{-\frac{n-\gamma}{2}} |u|^2 \big)u,
\end{equation*}
For dimension $n\geq 5$,  the exponent $\gamma=4$ is the unique
exponent which is energy critical in the sense that the natural
scale transformation
\begin{equation*}
u_{\lambda}(t,x)=\lambda^{\frac{n-2}{2}}u(\lambda^2 t, \lambda x),
\end{equation*}
leaves the energy invariant, in other words, the energy $E(u)$ is a
dimensionless quantity.

The Cauchy problem of the Hartree equation has been intensively
studied ([4-10], [15, 16, 18, 19]. With regard to the global
well-posedness and scattering results, they all dealt with the
$\dot{H}^1$-subcritical case $\big(2<\gamma<\min(4, n)\big)$ in
the energy space or some weighted spaces. In \cite{MiXZ}, we
obtained the small data scattering result for the
$\dot{H}^1$-critical case in the energy space. For the large
initial data for the $\dot{H}^1$-critical case $\big(\gamma=4,
n\geq5\big)$ in the energy space , the argument in \cite{MiXZ} can
not yield the global well-posedness, even with the conservation of
the energy $(\ref{cl})$, because the time of existence given by
the local theory depends on the profile of the data as well as on
the energy.

Concerning the $\dot{H}^1$-subcritical case $\big(2<\gamma<\min(4,
n)\big)$, using the method of Morawetz and Strauss \cite{MoStr72},
J. Ginibre and G. Velo \cite{GiV00} developed the scattering
theory in the energy space, where they exploited the properties of
$\Delta$ and obtain the usual Morawetz estimate
\begin{equation*}
\aligned -\int^{t_2}_{t_1} \int \int \big|u(t,x)\big|^2
\frac{x}{|x|} \frac{1}{|x-y|^{\gamma}}\nabla |u(t,y)|^2dydxdt
\lesssim CE(u).
\endaligned
\end{equation*}
Later, K. Nakanishi \cite{Na99d} exploited the properties of
$i\partial_t+\Delta$ and used a certain related Sobolev-type
inequality to obtain a new Morawetz estimate
\begin{equation*}
\aligned \int\int_{\mathbb{R}^{n+1}}
\frac{|t|^{1+\nu}|u(t,x)|^{\frac{2n}{n-2}}}{(|t|+|x|)^{2+\nu}}
dxdt \leq C(E, \nu), \quad \text{for any}\quad \nu>0,
\endaligned
\end{equation*}
 which was independent of the nonlinearity.

In this paper, we deal with the Cauchy problem of the Hartree
equation with the large data for the $\dot{H}^1$-critical case
$\big(\gamma=4, n\geq 5\big)$. Inspired by the approach of
Bourgain \cite{Bo98} and Tao \cite{Tao05} in the case of the
$\dot{H}^1$-critical Schr\"{o}dinger equation with the local
nonlinear term, we obtain the global well-posedness and scattering
results for the Hartree equation for the large radial data in
$\dot{H}^1$. The new ingredient is that we take advantage of the
following localized estimate for the first time
\begin{equation*}
\aligned - \int_{I}\int_{|x|\leq A|I|^{1/2}}|u(t,x)|^{2}\Delta
\Big(\frac{1}{|x|}\Big)dxdt =(n-3) \int_{I}\int_{|x|\leq A
|I|^{1/2}}\frac{|u(t,x)|^{2}}{|x|^{3}}dxdt \leq A|I|^{1/2}C(E)
\endaligned
\end{equation*}
 to rule out the possibility of
energy concentration, instead of the classical Morawetz estimate
\begin{equation*}
\aligned -\int_{I}\int\int
\Big(\frac{x}{|x|}-\frac{y}{|y|}\Big)\nabla
V(x-y)|u(x)|^{2}|u(y)|^{2}dydxdt
 \lesssim C(E)
\endaligned
\end{equation*}
due to the nonlinear term .

Our main result is the following global well-posedness result in the
energy space.

\begin{theorem}\label{main}
Let $n\geq 5$, and $\varphi \in \dot{H}^1$ be radial. then there
exists a unique global solution $u \in C^0_t (\dot{H}^1_x)\bigcap
L^6_tL^{\frac{6n}{3n-8}}_x$ to
\begin{equation} \label{equm}
\left\{ \aligned
    (iu_t +  \Delta u)(t,x)  & = \big(u V*|u|^2\big)  (t,x), \quad  \text{in}\  \mathbb{R}^n \times \mathbb{R}, \\
     u(0)&=\varphi(x), \quad \text{in} \ \mathbb{R}^n.
\endaligned
\right.
\end{equation}
where $V(x)=|x|^{-4}$ and on each compact time interval $[t_-,
t_+]$, we have
\begin{equation}\label{ine}
\Big\|u\Big\|_{L^6_tL^{\frac{6n}{3n-8}}_x([t_-, t_+]\times
\mathbb{R}^n)} \leq C(\big\| \varphi \big\|_{\dot{H}^1}).
\end{equation}
\end{theorem}

As the right hand side of $(\ref{ine})$ is independent of $t_-,
t_+$, we can obtain the global spacetime estimate. As a direct
consequence of the global $L^6_tL^{\frac{6n}{3n-8}}_x$ estimate,
we have scattering, asymptotic completeness, and uniform
regularity.

\begin{corollary}
Let $\varphi$ be radial and have finite energy. Then there exists
finite energy solutions $u_{\pm}(t,x)$ to the free Schr\"{o}dinger
equation $iu_t + \Delta u =0$ such that
\begin{equation*}
\big\|u_{\pm}(t)-u(t)\big\|_{\dot{H}^1}\rightarrow 0 \quad
\text{as}\quad t\rightarrow\pm \infty.
\end{equation*}
Furthermore, the maps $\varphi\mapsto u_{\pm}(0)$ are homeomorphisms
from $\dot{H}^1(\mathbb{R}^n)$ to $\dot{H}^1(\mathbb{R}^n)$.
Finally, if $\varphi\in H^s$ for some $s>1$, then $u(t)\in H^s$ for
all time $t$, and one has the uniform bounds
\begin{equation*}
\aligned
 \sup_{t\in \mathbb{R}} \big\|u(t)\big\|_{H^s} \leq C(E(\varphi),
 s)\big\|\varphi\big\|_{H^s}.
 \endaligned
\end{equation*}
\end{corollary}

The paper is organized as follows.

In Section $2$, we introduce notations and the basic estimates; In
Section $3$, we derive the local mass conservation and Morawetz
inequality; In Section $4$, we discuss the local theory for
$(\ref{equm})$; In Section $5$, we obtain the perturbation theory;
Finally, we prove the main theorem in Section $6$.

\section{Notations and basic estimates}
\setcounter{section}{2}\setcounter{equation}{0} We will often use
the notations $a\lesssim b$ and $a=O(b)$ to denote the estimate
$a\leq Cb$ for some $C$. The derivative operator $\nabla$ refers
to the space variable only. We also occasionally use subscripts to
denote the spatial derivatives and use the summation convention
over repeated indices.

We define $\langle a, b\rangle = \text{Re} (a\overline{b})$,
$\partial=(\partial_t, \nabla)$, $\mathcal{D}=(-\frac{i}{2},
\nabla)$; For $1\leq p\leq \infty,$ we denote by $p'$  the dual
exponent, that is, $\frac1p+ \frac1{p'} =1$.

For any time interval $I$, we use $L^q_tL^r_x(I\times \mathbb{R}^n)$
to denote the mixed spacetime Lebesgue norm
\begin{equation*}
\big\|u\big\|_{L^q_tL^r_x(I\times \mathbb{R}^n)}:=\bigg(
\int_I\big\| u\big\|^q_{L^r(\mathbb{R}^n)}dt \bigg)^{1/q}
\end{equation*}
with the usual modifications when $q=\infty$. When $q=r$, we
abbreviate $L^q_tL^r_x$ by $L^q_{t,x}$.

We use $U(t)=e^{it\Delta}$ to denote the free group generated by
the free Schr\"{o}dinger equation $iu_t+\Delta u=0$. It can
commute with derivatives, and obeys the inequality
\begin{equation}\label{dise}
\big\| e^{it\Delta}f \big\|_{L^p(\mathbb{R}^n)} \lesssim
|t|^{-n(\frac12-\frac1p)} \big\| f \big\|_{L^{p'}(\mathbb{R}^n)}
\end{equation}
for $t\not=0$, $2\leq p \leq \infty$.

We say that a pair $(q, r)$ is admissible if
\begin{equation*}
\frac{2}{q} = n\Big(\frac{1}{2}-\frac{1}{r}\Big),
\end{equation*}
and
\begin{equation*}
2\leq r \left\{
\begin{array}{rrl}
\leq & \infty, & n=1;\\
< &  \infty, & n=2; \\
\leq &\frac{2n}{n-2}, & n\geq 3.
\end{array}
\right.
\end{equation*}

For a spacetime slab $I\times \mathbb{R}^n$, we define the {\it
Strichartz} norm $\dot{S}^0(I)$ by
\begin{equation*}
\big\|u\big\|_{\dot{S}^0(I)}:= \sup_{(q, r)\ \text{admissible}}
\big\|u\big\|_{L^q_tL^r_x(I\times \mathbb{R}^n)}.
\end{equation*}
and define $\dot{S}^1(I)$ by
\begin{equation*}
\big\|u\big\|_{\dot{S}^1(I)}:= \big\|\nabla u\big\|_{\dot{S}^0(I)}.
\end{equation*}

When $n\geq 3$, the spaces $\big(\dot{S}^0(I),
\|\cdot\|_{\dot{S}^0(I)}\big)$ and $\big(\dot{S}^1(I),
\|\cdot\|_{\dot{S}^1(I)}\big)$ are Banach spaces, respectively.

Based on the above notations, we have the following {\it Strichartz
} inequalities

\begin{lemma}\cite{KeT98}, \cite{Stri77}\label{lemstri1}
Let $u$ be an $\dot{S}^0$ solution to the Schr\"{o}dinger equation
$(\ref{equ1})$. Then
\begin{equation*}
\big\|u\big\|_{\dot{S}^0} \lesssim \big\|u(t_0)\big\|_{L^2_x} +
\big\|f(u)\big\|_{L^{q'}_tL^{r'}_x(I\times \mathbb{R}^n)}
\end{equation*}
for any $t_0 \in I$ and any admissible pairs $(q, r)$. The implicit
constant is independent of the choice of interval $I$.
\end{lemma}

From Sobolev embedding, we have

\begin{lemma}\label{lemstri2}
For any function $u$ on $I\times \mathbb{R}^n$, we have
\begin{equation*}
\big\|\nabla u\big\|_{L^{\infty}_tL^2_x} + \big\|\nabla
u\big\|_{L^3_tL^{\frac{6n}{3n-4}}_x} + \big\|
u\big\|_{L^{\infty}_tL^{\frac{2n}{n-2}}_x} + \big\|
u\big\|_{L^6_tL^{\frac{6n}{3n-8}}_x} \lesssim \big\|
u\big\|_{\dot{S}^1},
\end{equation*}
where all spacetime norms are on $I\times \mathbb{R}^n$.
\end{lemma}

For convenience, we introduce two abbreviated notations. For a time
interval $I$, we denote
\begin{equation*}
\aligned
 \big\|u\big\|_{X(I)}  :=
\big\|u\big\|_{L^6_tL^{\frac{6n}{3n-8}}_x(I\times
\mathbb{R}^n)};\quad
 \big\|u\big\|_{W(I)} :=
\big\|\nabla u\big\|_{L^3_tL^{\frac{6n}{3n-4}}_x(I\times
\mathbb{R}^n)}.
\endaligned
\end{equation*}

\begin{lemma}\label{nle}
Let $\displaystyle f(u)(t,x)=\big(uV*|u|^2\big)(t,x)$, where
$V(x)=|x|^{-4}$. For any time interval $I$ and $t_0 \in I$, we have
\begin{equation*}
\Big\|\int^t_{t_0} e^{i(t-s)\Delta}f(u)(s,x) ds\Big\|_{\dot{S}^1(I)}
\lesssim \big\|u\big\|^2_{X(I)} \big\|u\big\|_{W(I)}.
\end{equation*}
\end{lemma}
{\bf Proof: } By Strichartz estimates, Hardy-Littlewood-Sobolev
inequality and H\"{o}lder inequality, we have
\begin{equation*}
\aligned & \quad \ \Big\|\int^t_{t_0} e^{i(t-s)\Delta}f(u)(s,x)
ds\Big\|_{\dot{S}^1(I)} \\
& \lesssim \|\nabla f(u)(t,x)\|_{L^\frac{3}{2}_tL^{\frac{6n}{3n+4}}_x(I\times \mathbb{R}^n)}\\
& \lesssim \|\nabla u
V*|u|^{2}\|_{L^\frac{3}{2}_tL^{\frac{6n}{3n+4}}_x(I\times
\mathbb{R}^n)}+\|u V*(u\nabla
u)\|_{L^\frac{3}{2}_tL^{\frac{6n}{3n+4}}_x(I\times \mathbb{R}^n)}\\
& \lesssim  \|\nabla u\|_{L^3_tL^{\frac{6n}{3n-4}}_x(I\times
\mathbb{R}^n)}\|V*|u|^{2}\|_{L^3_tL^{\frac{3n}{4}}_x(I\times
\mathbb{R}^n)} + \|u\|_{L^6_tL^{\frac{6n}{3n-8}}_x(I\times
\mathbb{R}^n)}\|V*(u\nabla
u)\|_{L^2_tL^{\frac{n}{2}}_x(I\times \mathbb{R}^n)} \\
&\lesssim \big\|u\big\|^2_{X(I)} \big\|u\big\|_{W(I)}.
\endaligned
\end{equation*}

\section{Local mass conservation and Morawetz inequality}
\setcounter{section}{3}\setcounter{equation}{0} In this section,
we will prove two useful estimates. One is a local mass
conservation estimate and the other is a Morawetz inequality,
which appears in Morawetz identity. The local mass conservation
estimate is used to control the flow of mass through a region of
space, and the Morawetz inequality is used to prevent
concentration.

\subsection{Local mass conservation}
We recall a local mass conservation law that has appeared in
\cite{Bo98}, \cite{KiVZ} and \cite{Tao05}. For completeness, we give
the sketch of the proof. Let $\chi$ be a bump function supported on
the ball $B(0, 1)$ that equals $1$ on the ball $B(0,1/2)$. Observe
that if $u$ is a finite energy solution of $(\ref{equm})$, then
\begin{equation*}
\partial_t \big| u(t,x)\big|^2 = -2 \nabla \cdot
\text{Im}(\overline{u}\nabla u(t,x)).
\end{equation*}

We define
\begin{equation*}
\aligned
 \text{Mass}(u(t), B(x_0, R)):= \int \Big|\chi\Big(\frac{x-x_0}{R}\Big)u(t,x)
 \Big|^2dx.
 \endaligned
\end{equation*}
Differentiating the above quantity with respect to time, we obtain
by the integration by parts
\begin{equation*}
\aligned
\partial_t \text{Mass}(u(t), B(x_0, R)) & = \int \Big|\chi\Big(\frac{x-x_0}{R}\Big)
 \Big|^2  \partial_t \big| u(t,x)\big|^2dx\\
 & = -2 \int \Big|\chi\Big(\frac{x-x_0}{R}\Big)
 \Big|^2  \nabla  \cdot
\text{Im}(\overline{u}\nabla u) dx\\
& = -\frac4R\ \int \chi\Big(\frac{x-x_0}{R}\Big) \nabla
\chi\Big(\frac{x-x_0}{R}\Big) \text{Im} (\overline{u}\nabla u) dx\\
& \lesssim \frac1R \big\|\nabla u(t)\big\|_{L^2} \Big(
\text{Mass}(u(t), B(x_0, R)) \Big)^{1/2},
\endaligned
\end{equation*}
hence, we have
\begin{equation}\label{meR}
\Big| \text{Mass}(u(t_1), B(x_0, R))^{1/2}  - \text{Mass}(u(t_2),
B(x_0, R))^{1/2} \Big| \lesssim \frac1R \big|t_1 -t_2\big|.
\end{equation}
This implies that if the local mass $\text{Mass}(u(t), B(x_0, R)) $
is large for some time $t$, then it can also be shown to be
similarly large for nearly time $t$, by increasing the radius $R$ if
necessary to reduce the rate of change of the mass.

On the other hand, from Sobolev and H\"{o}lder inequalities, we have
\begin{equation}\label{mer}
\aligned \text{Mass}(u(t), B(x_0, R))  \leq \Big\|
\chi\Big(\frac{x-x_0}{R}\Big)\Big\|^2_{L^n_x} \big\|u
\big\|^2_{L^{\frac{2n}{n-2}}_{x}} \lesssim R^2 \big\| \nabla
u\big\|^2_{L^2_x}.
\endaligned
\end{equation}
This gives the control of mass in small volumes.

\subsection{A Morawetz inequality}

To prevent the concentration of the energy, we need a Morawetz
estimate. The Morawetz estimate is based on some integral identity
derived by variation of the lagrangian.

We define $\ell(u)$ by
\begin{equation*}
2\ell(u)=\langle iu_t, u\rangle + |\nabla u|^2 + \frac12
|u|^2(V*|u|^2)
\end{equation*}
$\ell(u)$ is the lagrangian density associated to the equation
$(\ref{equ1})$.

From the definition of the variation of the functional $\ell$, we
have
\begin{equation*}
\aligned
\delta_{v}\ell(u):& =\lim_{\epsilon \rightarrow 0} \frac{\ell(u+\epsilon v)-\ell(u)}{\epsilon}\\
& = \langle iu_t + \Delta u -u(V*|u|^2) , v \rangle + \partial \cdot
\langle \mathcal{D}u, v \rangle.
\endaligned
\end{equation*}
Using this identity together with $\displaystyle h=\frac{x}{|x|},
q=\frac12 \text{Re}(\mathcal{D} \cdot h)=\frac{n-1}{2|x|}$ and
$Mu=h\cdot \mathcal{D}u+qu$, we obtain the following formula:
\begin{equation*}
\aligned \langle iu_t + \Delta u -u(V*|u|^2), Mu\rangle  = & -
\partial\cdot \langle \mathcal{D}u, Mu \rangle + \mathcal{D} \cdot \big( h\ell(u) + \frac{|u|^2}{2} \partial q
\big) \\
& + \sum^{n}_{\alpha =0} \langle \mathcal{D}u, \partial h_{\alpha}
\mathcal{D}_{\alpha}u\rangle - \frac{|u|^2}{2} \mathcal{D} \cdot
\partial q - (V* \nabla |u|^2)\cdot  \frac{x}{2|x|} |u|^2.
\endaligned
\end{equation*}

As a consequence of the above dilation identity, we have the
following Morawetz estimate, which plays an important role in our
proof.

\begin{proposition}[Morawetz estimate]\label{mor}
Let u be a solution to $(\ref{equm})$ on a spacetime slab $I\times
\mathbb{R}^n$. Then for any $A\geq 1$, we have
\begin{equation*}
\aligned &\int_{I}\int_{|x|\leq
A|I|^{1/2}}\frac{|u|^{2}}{|x|^{3}}dxdt-\int_{I}\int\int_{\Omega}
\Big(\frac{x}{|x|}-\frac{y}{|y|}\Big)\nabla
V(x-y)|u(x)|^{2}|u(y)|^{2}dydxdt\\
& \lesssim A|I|^{1/2}E,
\endaligned
\end{equation*}
where $\Omega=\big\{ (x, y)\in \mathbb{R}^n\times \mathbb{R}^n;
|x|\leq A|I|^{1/2}; |y|\leq A|I|^{1/2} \big\}.$
\end{proposition}

\begin{remark}
Since
\begin{equation*}
 -\Big(\frac{x}{|x|}-\frac{y}{|y|}\Big)\nabla
V(x-y)=4\frac{|x||y|-x\cdot y}{|x-y|^6}\Big( \frac{1}{|x|}+
\frac{1}{|y|}\Big)\geq 0,
\end{equation*}
we have
\begin{equation*}
-\int_{I}\int\int_{\Omega}
\Big(\frac{x}{|x|}-\frac{y}{|y|}\Big)\nabla
V(x-y)|u(x)|^{2}|u(y)|^{2}dydxdt \geq 0.
\end{equation*}
\end{remark}

{\bf Proof:} We define $V^a_0(t)=\displaystyle \int
a(x)|u(t,x)|^{2}dx$, then
$$M_{0}^{a}(t)=:\partial_{t}V^a_0(t)=2\mathrm{Im}\int a_{j}u_{j}\overline{u}dx.$$
and
\begin{equation*}
\aligned
\partial_{t}M_{0}^{a}(t)&=-2\mathrm{Im}\int
a_{jj}u_{t}\overline{u}dx-4\mathrm{Im}\int a_{j}\overline{u}_{j}u_{t}dx\\
&=-\int\triangle\triangle a|u|^{2}dx+4\mathrm{Re}\int
a_{jk}\overline{u_{j}}u_{k}dx\\
& \qquad -2\mathrm{Re}\int\int \nabla a(x)\nabla
V(x-y)|u(y)|^2|u(x)|^{2}dxdy\\
&=-\int\triangle\triangle a|u|^{2}dx+4\mathrm{Re}\int
a_{jk}\overline{u_{j}}u_{k}dx \\
& \qquad -\mathrm{Re}\int\int \big(\nabla a(x)- \nabla
a(y)\big)\nabla
V(x-y)|u(y)|^2|u(x)|^{2}dxdy\\
\endaligned
\end{equation*}
where we use the symmetry of $a(x)$ and $V(x)$. Let $R>0$ and let
$\eta$ be a bump function adapted to the ball $|x|\leq R$ which
equals 1 on the ball $|x|\leq R/2$. We set $a(x):=|x|\eta(x)$.

For $|x|\leq R/2$, we have
$$a_j=\frac{x_{j}}{|x|};\quad a_{jk}=\frac{\delta_{jk}}{|x|}-\frac{x_{j}x_{k}}{|x|^{3}};\quad
\triangle a=\frac{n-1}{|x|};\quad -\triangle\triangle
a=\frac{(n-1)(n-3)}{|x|^{3}}.$$ and for $R/2\leq |x|\leq R$, we have
bounds
$$a_{j}=O(1);\quad a_{jk}=O(R^{-1});\quad \triangle\triangle
a=O(R^{-3}).$$ Thus we have
\begin{equation*}
\aligned
\partial_{t}M_{0}^a(t)&=(n-1)(n-3)\int_{|x|\leq
R/2}\frac{|u|^{2}}{|x|^{3}}dx+4\int_{|x|\leq R/2}\frac{|\nabla u|^{2}- |\partial_r u|^2}{|x|}dx\\
&-\int\int_{\Omega_1} \Big(\frac{x}{|x|}-\frac{y}{|y|}\Big)\nabla
V(x-y)|u(x)|^{2}|u(y)|^{2}dydx \\
&+O\Big(\int_{|x|\sim R}\Big(\frac{|u|^{2}}{R^{3}}+\frac{|\nabla
u|^{2}}{R}\Big)dx\Big)\\
& +O\Big(\int\int_{\Omega_2} \big(a_j(x)-a_j(y) \big)
\frac{x_j-y_j}{|x-y|^{\gamma+2}} |u(x)|^2 |u(y)|^2dydx\Big)
\endaligned
\end{equation*}
where $\gamma =4$,
\begin{equation*}
\aligned \Omega_1 & =\big\{ (x,y)\in \mathbb{R}^n\times
\mathbb{R}^n; |x|\leq R/2, |y|\leq R/2 \big\};\\
\Omega_2&=\big\{(x,y)\in \mathbb{R}^n\times \mathbb{R}^n;
|x|\thicksim R \big\} \cup \big\{(x,y)\in \mathbb{R}^n\times
\mathbb{R}^n; |y|\thicksim R \big\}.\\
\endaligned
\end{equation*}
Meanwhile
$$\int_{|x|\sim R}\Big(\frac{|u|^{2}}{R^{3}}+\frac{|\nabla
u|^{2}}{R}\Big)dx\lesssim R^{-1}E,$$
\begin{equation*}
\aligned &\quad\  \Big|\int\int_{\Omega_2} \big(a_j(x)-a_j(y)
\big)
\frac{x_j-y_j}{|x-y|^{\gamma+2}} |u(x)|^2 |u(y)|^2dydx \Big| \\
&\leq \Big|\int\int_{\Omega_2: \ |x-y|\leq R/4} \big(a_j(x)-a_j(y)
\big)
\frac{x_j-y_j}{|x-y|^{\gamma+2}} |u(x)|^2 |u(y)|^2dydx \Big| \\
&\qquad \qquad + \Big|\int\int_{\Omega_2: \ |x-y|\geq R/4}
\big(a_j(x)-a_j(y) \big) \frac{x_j-y_j}{|x-y|^{\gamma+2}} |u(x)|^2
|u(y)|^2dydx \Big|\\
& \lesssim R^{-1} \Big|\int\int_{\Omega_2}
\frac{1}{|x-y|^{\gamma}} |u(x)|^2 |u(y)|^2dydx \Big| \\
 &
\lesssim R^{-1}E.
\endaligned
\end{equation*} Moreover, from Sobolev and H\"{o}lder
inequalities, we have
$$M_0^a(t)\lesssim\int_{|x|\lesssim R}|u||\nabla u|\lesssim
\|u\|_{L^{\frac{2n}{n-2}}_x}\|\nabla
u\|_{L^2_x}\Big(\int_{|x|\lesssim R}dx\Big)^{1/n}\lesssim RE.$$ So
if we integrate by parts on a time interval I and take $R=2
A|I|^{1/2}$, we obtain
\begin{equation*}
\aligned &\int_{I}\int_{|x|\leq
A|I|^{1/2}}\frac{|u|^{2}}{|x|^{3}}dxdt-\int_{I}\int\int_{\Omega}
\Big(\frac{x}{|x|}-\frac{y}{|y|}\Big)\nabla
V(x-y)|u(x)|^{2}|u(y)|^{2}dydxdt\\
& \lesssim A|I|^{1/2}E
\endaligned
\end{equation*}
for $n\geq 4$. The proof is completed.

\section{Local theory}
\setcounter{section}{4}\setcounter{equation}{0} In this section, we
develop a local well-posedness and blow-up criterion for the
$\dot{H}^1$-critical Hartree equation. First, we have

\begin{proposition}[Local well-posedness]\label{lwp}
Let $u(t_0) \in \dot{H}^1$, and $I$ be a compact time interval that
contains $t_0$ such that
\begin{equation*}
\big\| U(t-t_0) u(t_0) \big\|_{X(I)}\leq \eta,
\end{equation*}
for a sufficiently small absolute constant $\eta >0$. Then there
exists a unique strong solution to $(\ref{equm})$ on $I\times
\mathbb{R}^n$ such that
\begin{equation*}
\big\| u\big\|_{X(I)} \leq C(\big\|u(t_0)\big\|_{\dot{H}^1}).
\end{equation*}
\end{proposition}
{\bf Proof: } The proof of this proposition is standard and based on
the contraction mapping arguments. We define the solution map to
be$$\Phi(u)(t):=U(t-t_0)u(t_0)-i\int_{t_0}^tU(t-s)f(u(s))ds,$$ then
$\Phi$ is a map from
\begin{equation*}\mathcal{B}=\{u:\|u\|_{X(I)}\leq2\eta,
\|u\|_{W(I)}\leq2 C \|u(t_0)\|_{\dot{H}^{1}}\}
\end{equation*}
with the metric
$$\|u\|_{\mathcal{B}}=\|u\|_{X(I)}+\|u\|_{W(I)}$$onto itself
because $$\|\Phi(u)\|_{X(I)}\leq
\|U(t-t_0)u(t_0)\|_{X(I)}+C\|u\|_{X(I)}^2\|u\|_{W(I)}\leq
\eta+8C\eta^2\|u(t_0)\|_{\dot{H}^{1}}\leq2\eta;$$
$$\|\Phi(u)\|_{W(I)}\leq C \|u(t_0)\|_{\dot{H}^{1}}+C\|u\|_{X(I)}^2
\|u\|_{W(I)}\leq
C\|u(t_0)\|_{\dot{H}^{1}}+8C\eta^2\|u(t_0)\|_{\dot{H}^{1}}\leq2C\|u(t_0)\|_{\dot{H}^{1}}.$$
 It
suffices to prove $\Phi$ is a contraction map. Let $u$, $v\in
\mathcal{B}$, then
\begin{equation*}
\aligned
\|\Phi(u)-\Phi(v)\|_{W(I)}\leq& \Big\|\int_{0}^{t}U(t-s)(V*(\bar{u}-\bar{v}) u)u(s,x)ds\Big\|_{W(I)}\\
&+\Big\|\int_{0}^{t}U(t-s)(V*\bar{v}( u-
v))u(s,x)ds\Big\|_{W(I)}\\
&+\Big\|\int_{0}^{t}U(t-s)(V*\bar{v} v)(u-v)(s,x)ds\Big\|_{W(I)}.
\endaligned
\end{equation*}
By Lemma \ref{nle}, we have
\begin{equation*}
\aligned \|\Phi(u)-\Phi(v)\|_{W(I)}&\leq
\|u-v\|_{X(I)}\bigg(\|u\|_{W(I)}\|u\|_{X(I)}+\|u\|_{W(I)}\|v\|_{X(I)}+ \|v\|_{W(I)}\|v\|_{X(I)}\bigg)\\
& \qquad +
\|u-v\|_{W(I)}\bigg(\|u\|_{X(I)}\|u\|_{X(I)}+\|u\|_{X(I)}\|v\|_{X(I)}+ \|v\|_{X(I)}\|v\|_{X(I)}\bigg)\\
&\leq 12C\eta\|u(t_0)\|_{\dot{H}^{1}}\|u-v\|_{X(I)}+12\eta^2\|u-v\|_{W(I)}\\
&\leq\frac{1}{4}(\|u-v\|_{X(I)}+\|u-v\|_{W(I)})
\endaligned
\end{equation*}
In the same way, we have
\begin{equation*}
\aligned
\|\Phi(u)-\Phi(v)\|_{X(I)}&\leq 12C\eta\|u(t_0)\|_{\dot{H}^{1}}\|u-v\|_{X(I)}+12\eta^2\|u-v\|_{W(I)}\\
&\leq\frac{1}{4}(\|u-v\|_{X(I)}+\|u-v\|_{W(I)})
\endaligned
\end{equation*}
as long as $\eta$ is chosen sufficiently small. Then the contraction
mapping theorem implies the existence of the unique solution to
(\ref{equm}) on I.

Next, we give the blow-up criterion of the solutions for
$(\ref{equm})$. The usual form is similar to those in \cite{Ca03},
\cite{KeM06s}, which is in the form of a maximal interval of
existence. For convenience, we obtain

\begin{proposition}[Blow-up criterion]\label{blowup}
Let $\varphi \in \dot{H}^1$, and let $u$ be a strong solution to
$(\ref{equm})$ on the slab $[0, T)\times \mathbb{R}^n$ such that
\begin{equation*}
\big\| u \big\|_{X([0,T))}< \infty.
\end{equation*}
Then there exists $\delta>0$ such that the solution $u$ extends to a
strong solution to $(\ref{equm})$ on the slab $[0, T+\delta]\times
\mathbb{R}^n$.
\end{proposition}

{\bf Proof: } By the absolute continuity of integrals, there exists
a $t_{0}\in [0,T)$, such that
$$\|u\|_{X([t_{0}, T))}\leq\eta/4,$$ then by Lemma \ref{nle}, we
have
\begin{equation*}\|u\|_{W([t_{0},
T))}\lesssim\|u(t_0)\|_{\dot{H}^{1}}+\|u\|_{X([t_{0},
T))}^{2}\|u\|_{W([t_{0}, T))},\end{equation*}therefore
\begin{equation*}\|u\|_{W([t_{0}, T))}\lesssim\|u(t_0)\|_{\dot{H}^{1}}.
\end{equation*}

Now we write
\begin{equation*}U(t-t_0)u(t_0)=u(t)+i\int_{t_0}^{t}U(t-s)(V*|u|^2)u(s,x)ds,
\end{equation*}
then \begin{equation*}\|U(t-t_0)u(t_0)\|_{X([t_{0}, T))}\leq
\|u\|_{X([t_{0}, T))}+C\|u\|_{X([t_{0}, T))}^{2}\|u\|_{W([t_{0},
T))}\leq\frac{\eta}{4}+C\eta^2\|u(t_0)\|_{\dot{H}^{1}}\leq\frac{\eta}{2}.
\end{equation*}
By the absolute continuity of integrals again, there exists a
$\delta$, such that
$$\|U(t-t_0)u(t_0)\|_{X( [t_0, T+\delta))}\leq\eta.$$
Thus we may apply Proposition \ref{lwp} on the interval $[t_0,
T+\delta]$ to complete the proof.

In other words, this lemma asserts that if $[t_0,T^*)$ is the
maximal interval of existence and $T^*<\infty$, then
$$\|u\|_{X([t_0,T^*))}=\infty.$$

\section{Perturbation result}
\setcounter{section}{5}\setcounter{equation}{0}

In this section, we obtain the perturbation for Hartree equation,
which shows that the solution can not be large if the linear part
of the solution is not large. This is an analogue of Lemma $3.2$
in \cite{Tao05}, and later, Killip, Visan and Zhang \cite{KiVZ}
gave the similar perturbation result for the Schr\"{o}dinger
equation with the quadric potentials.
\begin{lemma}[Perturbation lemma]\label{pertl}
Let $u$ be a solution to $(\ref{equm})$ on $I=[t_1, t_2]$ such that
\begin{equation*}
\frac12 \eta \leq \big\|u\big\|_{X(I)} \leq \eta,
\end{equation*}
where $\eta$ is sufficiently small constant depending on the norm of
the initial data, then
\begin{equation*}
\aligned \big\| u\big\|_{\dot{S}^1(I)} \lesssim 1,\quad
\big\|u_k\big\|_{X(I)}  \geq \frac14 \eta,
 \endaligned
\end{equation*}
where $u_k(t)=U(t-t_k)u(t_k)$ for $k=1,2$.
\end{lemma}

{\bf Proof: } From Strichartz estimate and Lemma $\ref{nle}$, we
obtain
\begin{equation*}
\aligned \big\|u\big\|_{\dot{S}^1(I)} & \lesssim
\big\|u(t_1)\big\|_{\dot{H}^1} + \big\|u\big\|^2_{X(I)}
\big\|u\big\|_{W(I)} \\
& \lesssim \big\|u(t_1)\big\|_{\dot{H}^1} + \big\|u\big\|^2_{X(I)}
\big\|u\big\|_{\dot{S}^1(I)} \\
&  \lesssim \big\|u(t_1)\big\|_{\dot{H}^1} + \eta^2
\big\|u\big\|_{\dot{S}^1(I)}. \\
\endaligned
\end{equation*}

If $\eta$ is sufficiently small, we have the first claim
\begin{equation*}
\aligned \big\| u\big\|_{\dot{S}^1(I)}  \lesssim 1.
 \endaligned
\end{equation*}

As for the second claim, we give the proof for $k=1$, the case $k=2$
is similar. Using Strichartz estimate and Lemma $\ref{nle}$ again,
we have
\begin{equation*}
\aligned \big\|u-u_1\big\|_{X(I)}  \lesssim \eta^2
\big\|u\big\|_{\dot{S}^1(I)} \lesssim  \eta^2,
\endaligned
\end{equation*}
therefore, the second claim follows by the triangle inequality and
choosing $\eta$ sufficiently small.

\section{Global well-posedness}
\setcounter{section}{6}\setcounter{equation}{0} In this section,
we give the proof of Theorem $\ref{main}$. The new ingredient is
that we first take advantage of the the estimate of the term
$\displaystyle \int_{I}\int_{|x|\leq
A|I|^{1/2}}\frac{|u|^{2}}{|x|^{3}}dxdt$ in the localized Morawetz
identity to rule out the possibility of energy concentration,
which is independent of the nonlinear term. For the
Schr\"{o}dinger equation, Tao \cite{Tao05} used the classical
Morawetz estimate, which depends on the nonlinearity,  to prevent
the concentration.

For readability, we first take some constants
\begin{equation}\label{constant}
C_1=6n;  \quad C_2=3 ;\quad C_3=18n.
\end{equation}
which come from several constraints in the rest of this section.
All implicit constants in this section are permitted to depend on
the dimension $n$ and the energy.

Fix $E$, $[t_-, t_{+}]$, $u$. We may assume that the energy is
large, $E>c>0$, otherwise the claim follows from the small energy
theory \cite{MiXZ}. From the boundedness of energy and Sobolev
embedding, we can obtain
\begin{equation}\label{eb}
\big\|u(t)\big\|_{\dot{H}^1_x} +
\big\|u(t)\big\|_{L^{\frac{2n}{n-2}}_x} \lesssim 1
\end{equation}
for all $t \in [t_-, t_{+}]$.

Assume that the solution $u$ already exists on $[t_-, t_{+}]$. By
Lemma $\ref{blowup}$, it suffices to obtain a priori estimate
\begin{equation}\label{object1}
\big\|u\big\|_{X([t_-, t_+])} \leq O(1),
\end{equation}
where $O(1)$ is independent of $t_-$, $t_+$.

We may assume that
\begin{equation*}
\big\|u\big\|_{X([t_-, t_+])} \geq 2\eta,
\end{equation*}
otherwise it is trivial. We divide $[t_-, t_+]$ into $J$
subintervals $I_j=[t_j, t_{j+1}]$ for some $J\geq 2$ such that
\begin{equation}\label{ij}
\frac\eta2 \leq \big\|u\big\|_{X(I_j)}\leq \eta,
\end{equation}
where $\eta$ is a small constant depending on the dimension $n$ and
the energy. As a consequence, it suffices to estimate the number
$J$.

Now let $u_{\pm} = U(t-t_\pm)u(t_{\pm})$. By Sobolev embedding and
Strichartz estimates, we have
\begin{equation}\label{freee}
\big\|u_{\pm}\big\|_{X([t_-, t_+])} \lesssim 1.
\end{equation}

We adapt the following definition of Tao \cite{Tao05}.
\begin{definition}
We call $I_j$ exceptional if
\begin{equation*}
\big\|u_{\pm}\big\|_{X(I_j)} > \eta^{C_3}
\end{equation*}
for at least one sign $\pm$. Otherwise, we call $I_j$ unexceptional.
\end{definition}

From $(\ref{freee})$, we obtain the upper bound on the number of
exceptional intervals, $O(\eta^{-6C_3})$. We may assume that there
exist unexceptional intervals, otherwise the claim would follow
from this bound and $(\ref{ij})$. Therefore, it suffices to
compute the number of unexceptional intervals.

We first prove the existence of a bubble of mass concentration in
each unexceptional interval.

\begin{proposition}[Existence of a bubble]\label{massconc1}
Let $I_j$ be an unexceptional interval. Then there exists $x_j\in
\mathbb{R}^n$ such that
\begin{equation*}
\mathrm{Mass}(u(t), B(x_j, \eta^{-C_1}|I_j|^{1/2})) \gtrsim
\eta^{C_1}|I_j|
\end{equation*}
for all $t\in I_j$.
\end{proposition}

{\bf Proof: } By time translation invariance and scale invariance,
we may assume that $I_j=[0,1]$. We subdivide $I_j$ further into $[0,
\frac12]$ and $[\frac12, 1]$. By $(\ref{ij})$ and the pigeonhole
principle and time reflection symmetry if necessary, we may assume
that
\begin{equation*}
\aligned \big\|u\big\|_{X([\frac12, 1])} \geq \frac{\eta}{4}.
\endaligned
\end{equation*}
Thus by Lemma $\ref{pertl}$, we have
\begin{equation}\label{pertla}
\Big\|U(t-\frac12)u(\frac12)\Big\|_{X([\frac12, 1])}\geq
\frac{\eta}{8}.
\end{equation}

By {\it Duhamel} formula, we have
\begin{equation}\label{duhamel}
\aligned
U(t-\frac12)u(\frac12)=U(t-t_-)u(t_-)&-i\int^\frac12_{0}U(t-s)f(u(s))ds\\
& -i\int^0_{t_-}U(t-s)f(u(s))ds.
\endaligned
\end{equation}

Since $[0, 1]$ is unexceptional interval, we have
\begin{equation*}
\big\| U(t-t_-)u(t_-) \big\|_{X([\frac12, 1])} = \big\| u_-(t)
\big\|_{X([\frac12, 1])} \leq \eta^{C_3}.
\end{equation*}

On the other hand, by $(\ref{ij})$, Lemma $\ref{lemstri2}$ , Lemma
$\ref{nle}$ and Lemma $\ref{pertl}$, we have
\begin{equation*}
\aligned \Big\|\int^\frac12_{0}U(t-s)f(u(s))ds \Big\|_{X([\frac12,
1])} & \lesssim \big\|u\big\|^2_{X([\frac12, 1])}
\big\|u\big\|_{W([\frac12, 1])}\\
& \lesssim \eta^2\big\|u\big\|_{\dot{S}^1([\frac12, 1])}  \lesssim
\eta^2.
\endaligned
\end{equation*}
Thus the triangle inequality implies that
\begin{equation*}
\aligned \Big\|\int^0_{t_-}U(t-s)f(u(s))ds \Big\|_{X([\frac12, 1])}
& \geq \frac{1}{100} \eta,
\endaligned
\end{equation*}
provided $\eta$ is chosen sufficiently small. Hence, if we define
\begin{equation*}
v(t):= \int^0_{t_-}U(t-s)f(u(s))ds,
\end{equation*}
then we have
\begin{equation}\label{ve}
\big\|v\big\| _{X([\frac12, 1])}  \geq \frac{1}{100} \eta.
\end{equation}

Next, we estimate the upper bound on $v$. We have by
$(\ref{duhamel})$ and the triangle inequality
\begin{equation}\label{ve2}
\aligned \big\|v\big\|_{\dot{S}^1([\frac12, 1])} & \leq
\Big\|U\big(t-\frac12\big)u(\frac12) \Big\|_{\dot{S}^1([\frac12,
1])} +
\big\|U(t-t_-)u(t_-) \big\|_{\dot{S}^1([\frac12, 1])} \\
& \qquad  + \Big\| \int^\frac12_{0}U(t-s)f(u(s))ds
\Big\|_{\dot{S}^1([\frac12, 1])} \\
& \lesssim \big\|u\big\|_{\dot{H}^1} + \big\|u\big\|^2_{X([0,
\frac12])} \big\|u\big\|_{W([0, \frac12])} \\
& \lesssim \big\|u\big\|_{\dot{H}^1} + \big\|u\big\|^2_{X([0,
\frac12])} \big\|u\big\|_{\dot{S}^1([0, \frac12])} \\
& \lesssim 1.
\endaligned
\end{equation}
where we use Strichartz estimate, $(\ref{ij})$ and Lemma
$\ref{pertl}$.

We shall need some additional regularity control on $v$. For any
$h\in \mathbb{R}^n$, let $u^{(h)}$ denote the translation of $u$ by
$h$, i.e. $u^{(h)}(t, x)=u(t, x-h)$.

\begin{lemma} \label{vreg}
Let $\chi$ be a bump function supported on the ball $B(0, 1)$ of
total mass one, and define
\begin{equation*}
\aligned
 v_{av}(t,x)=\int \chi(y)v(t,x+\eta^{C_2}y)dy,
\endaligned
\end{equation*}
then we have
\begin{equation*}
\big\|v-v_{av}\big\|_{X([\frac12, 1])} \lesssim \eta^{C_2}.
\end{equation*}
\end{lemma}

{\bf Proof: } By the chain rule, H\"{o}lder inequality and Sobolev
embedding, we have
\begin{equation*}
\aligned \big\| \nabla f(u)(s)\big\|_{L^{\frac{2n}{n+4}}_x} & \leq
\big\| (V*|u|^2) \nabla u\big\|_{L^{\frac{2n}{n+4}}_x} +
\big\| u(V*\nabla |u|^2)\big\|_{L^{\frac{2n}{n+4}}_x} \\
& \leq \big\| \nabla u \big\|_{L^2} \big\| V*
|u|^2\big\|_{L^{\frac{n}{2}}} + \big\|u\big\|_{L^{\frac{2n}{n-2}}}
\big\|V*\nabla |u|^2 \big\|_{L^{\frac{n}{3}}}\\
& \leq \big\| \nabla u \big\|_{L^2} \big\|
|u|^2\big\|_{L^{\frac{n}{n-2}}} + \big\|u\big\|_{L^{\frac{2n}{n-2}}}
\big\|\nabla |u|^2 \big\|_{L^{\frac{n}{n-1}}}\\
& \lesssim 1,
\endaligned
\end{equation*}
it follows by $(\ref{dise})$
\begin{equation*}
\aligned \big\| \nabla
v\big\|_{L^{\infty}_tL^{\frac{2n}{n-4}}_x([\frac12, 1]\times
\mathbb{R}^n)} \leq \sup_{t\in [\frac12, 1] }\int^0_{t_-}
\frac{1}{|t-s|^2} \big\| \nabla f(u)(s)\big\|_{L^{\frac{2n}{n+4}}_x}
ds \lesssim 1.
\endaligned
\end{equation*}
From $(\ref{ve2})$ and interpolation, we have
\begin{equation*}
\aligned \big\| \nabla
v\big\|_{L^{\infty}_tL^{\frac{6n}{3n-8}}_x([\frac12, 1]\times
\mathbb{R}^n)} & \leq \big\| \nabla
v\big\|^{2/3}_{L^{\infty}_tL^{\frac{2n}{n-4}}_x([\frac12, 1]\times
\mathbb{R}^n)} \big\| \nabla
v\big\|^{1/3}_{L^{\infty}_tL^2_x([\frac12, 1]\times \mathbb{R}^n)}
\\
& \lesssim 1.
\endaligned
\end{equation*}

From the fundamental theorem of calculus, we have
\begin{equation*}
\aligned \big\|
v-v^{(h)}\big\|_{L^{\infty}_tL^{\frac{6n}{3n-8}}_x([\frac12,
1]\times \mathbb{R}^n)} & \lesssim |h|.
\endaligned
\end{equation*}
This implies
\begin{equation*}
\aligned
\big\|v-v_{av}\big\|_{L^{\infty}_tL^{\frac{6n}{3n-8}}_x([\frac12,
1]\times \mathbb{R}^n)} & \leq \int \chi(y)\big\|
v(t,x+\eta^{C_2}y) -v(x)
\big\|_{L^{\infty}_tL^{\frac{6n}{3n-8}}_x([\frac12, 1]\times
\mathbb{R}^n)} dy \\
& \lesssim \int \chi(y) |\eta^{C_2}y|dy\\
& \lesssim \eta^{C_2}.
\endaligned
\end{equation*}
Hence from H\"{o}lder inequality, we obtain
\begin{equation*}
\aligned \big\|v-v_{av}\big\|_{X([\frac12, 1])}  \lesssim
\big\|v-v_{av}\big\|_{L^{\infty}_tL^{\frac{6n}{3n-8}}_x([\frac12,
1]\times \mathbb{R}^n)} \lesssim \eta^{C_2}.
\endaligned
\end{equation*}

This completes the proof of Lemma.

Now we return to the proof of Proposition $\ref{massconc1}$. By
Lemma $\ref{vreg}$ and $(\ref{ve})$, we have
\begin{equation}\label{vav}
\big\|v_{av}\big\|_{X([\frac12, 1])}\gtrsim \eta.
\end{equation}

On the other hand, by H\"{o}lder inequality, Young inequalities and
$(\ref{ve2})$, we have
\begin{equation*}
\aligned
\big\|v_{av}\big\|_{L^{\frac{2(3n-8)}{n-2}}_tL^{\frac{2n}{n-2}}_x([\frac12,
1]\times \mathbb{R}^n)} & \lesssim
\big\|v_{av}\big\|_{L^{\infty}_tL^{\frac{2n}{n-2}}_x([\frac12,
1]\times \mathbb{R}^n)} \\
&\lesssim \big\|v\big\|_{L^{\infty}_tL^{\frac{2n}{n-2}}_x([\frac12,
1]\times \mathbb{R}^n)} \\
&\lesssim  1.
\endaligned
\end{equation*}
Interpolating with $(\ref{vav})$ gives
\begin{equation*}
\aligned \big\|v_{av}\big\|_{L^{\infty}_{t,x}([\frac12, 1]\times
\mathbb{R}^n))}  & \gtrsim
\big\|v_{av}\big\|^{\frac{3n-6}{2}}_{X([\frac12, 1])}
\big\|v_{av}\big\|^{-\frac{3n-8}{2}}_{L^{\frac{2(3n-8)}{n-2}}_tL^{\frac{2n}{n-2}}_x([\frac12,
1]\times \mathbb{R}^n)}\\
& \gtrsim \eta^{\frac{3n-6}{2}}.
\endaligned
\end{equation*}
Thus there exists $(s_j, x_j)\in [\frac12, 1]\times \mathbb{R}^n$
such that
\begin{equation*}
\big|v_{av}(s_j, x_j)\big| \gtrsim \eta^{\frac{3n-6}{2}}.
\end{equation*}
Hence, by Cauchy-Schwarz inequality, we have
\begin{equation*}
\aligned  \big|v_{av}(s_j, x_j)\big| & = \Big|\int \chi(y)v(s_j,
x_j+\eta^{C_2}y)dy\Big|\\
& =\eta^{-nC_2}\Big|\int \chi(\frac{x-x_j}{\eta^{C_2}})v(s_j,x)dx\Big|\\
& \lesssim \eta^{-nC_2} \eta^{\frac{n}{2}C_2} \text{Mass}(v(s_j),
B(x_j, \eta^{C_2}))^{1/2},
\endaligned
\end{equation*}

that is
\begin{equation}\label{vcon1}
\text{Mass}(v(s_j), B(x_j, \eta^{C_2})) \gtrsim \eta^{3n-6+nC_2}
\gtrsim \eta^{C_1}.
\end{equation}

Observe that $(\ref{meR})$ also holds for $v$. If we take
$R=\eta^{-C_1}$ and choose $\eta$ sufficiently small, we have
\begin{equation}\label{vcon2}
\aligned
 \text{Mass}(v(t), B(x_j, \eta^{-C_1}))  & \gtrsim \big(
\text{Mass}(v(s_j),
B(x_j, \eta^{-C_1}))^{1/2} -\frac{1}{\eta^{-C_1}} \big)^2 \\
& \gtrsim ( \text{Mass}(v(s_j), B(x_j, \eta^{C_2}))^{1/2}
-\eta^{C_1}
)^2 \\
& \gtrsim\eta^{C_1}
\endaligned
\end{equation}
for all $t\in [0, 1]$.

The last step is to show that this mass concentration holds for $u$.
We first show mass concentration for $u$ at time $0$.

Since $[0, 1]$ is unexceptional interval, by the pigeonhole
principle, there is a $\tau_j \in [0, 1]$ such that
\begin{equation*}
\big\|u_-(\tau_j)\big\|_{L^{\frac{6n}{3n-8}}_x} \lesssim \eta^{C_3},
\end{equation*}
and so by H\"{o}lder inequality,
\begin{equation*}
\aligned
 \text{Mass}(u_-(\tau_j), B(x_j, \eta^{-C_1}))  & \lesssim
 \Big\|\chi\Big(\frac{x-x_j}{\eta^{-C_1}}\Big)\Big\|^2_{L^{\frac{3n}{4}}_x}
 \big\|u_-(\tau_j)\big\|^2_{L^{\frac{6n}{3n-8}}_x}\\
 & \lesssim \eta^{-\frac{8}{3}C_1+2C_3} \lesssim \eta^{2C_1}.
\endaligned
\end{equation*}

From $(\ref{meR})$, we have
\begin{equation}\label{ulcon}
\aligned
 \text{Mass}(u_-(0), B(x_j, \eta^{-C_1}))  \lesssim \eta^{2C_1}.
\endaligned
\end{equation}

Recall that $u(0)=u_-(0)-iv(0)$. Combing $(\ref{vcon2})$ and
$(\ref{ulcon})$ with the triangle inequality, we obtain

\begin{equation}\label{ucon3}
\aligned
 \text{Mass}(u(0), B(x_j, \eta^{-C_1}))  \gtrsim \eta^{C_1}.
\endaligned
\end{equation}
Using $(\ref{meR})$ again, we obtain the result.

Next, we use the radial assumption to show that the bubble of mass
concentration must occur at the spatial origin. In the forthcoming
paper, we shall use the interaction Morawetz estimate with the
frequency localized $L^2$ almost-conservation law to rule out the
possibility of the energy concentration at any place and deal with
the non-radial data. The corresponding results for the
Schr\"{o}dinger equation with local nonlinearity, please see
\cite{CKSTT07}, \cite{RyV05} and \cite{Vi05}.

\begin{corollary}[Bubble at the origin]\label{massconc2}
Let $I_j$ be an unexceptional interval. Then
\begin{equation*}
\mathrm{Mass}(u(t), B(0, \eta^{-3C_1}|I_j|^{1/2})) \gtrsim
\eta^{C_1}|I_j|
\end{equation*}
for all $t\in I_j$.
\end{corollary}

{\bf Proof: } If $x_j$ in Proposition $\ref{massconc1}$ is within
$\frac12 \eta^{-3C_1}|I_j|^{1/2}$ of the origin, then the result
follows immediately. Otherwise by the radial assumption, there would
be at least
\begin{equation*}
\aligned
O\Big(\frac{(\eta^{-3C_1}|I_j|^{1/2})^{n-1}}{(\eta^{-C_1}|I_j|^{1/2})^{n-1}}
\Big) \thickapprox O\big( \eta^{-2(n-1)C_1}\big)
\endaligned
\end{equation*}
many distinct balls each containing at least $\eta^{C_1} |I_j|$
amount of mass. By H\"{o}lder inequality, this implies
\begin{equation*}
\aligned \eta^{-2(n-1)C_1} \times \eta^{C_1} |I_j| & \lesssim
\int_{(\eta^{-3C_1}-\eta^{-C_1})|I_j|^{1/2}\leq |x| \leq
(\eta^{-3C_1}+\eta^{-C_1})|I_j|^{1/2}} |u(t,x)|^2dx \\
& \lesssim
\big\|u\big\|^2_{L^{\frac{2n}{n-2}}_x} \times \Big(\int_{(\eta^{-3C_1}-\eta^{-C_1})|I_j|^{1/2}\leq |x| \leq (\eta^{-3C_1}+\eta^{-C_1})|I_j|^{1/2}}dx \Big)^{2/n}\\
&\thickapprox \big\|u\big\|^2_{L^{\frac{2n}{n-2}}_x} \times
\Big(\big( \eta^{-3C_1} |I_j|^{1/2}\big)^{n-1} \times \eta^{-C_1}
|I_j|^{1/2}\Big)^{\frac2n},
\endaligned
\end{equation*}
that is
\begin{equation*}
\aligned \big\|u\big\|^2_{L^{\frac{2n}{n-2}}_x}  \gtrsim
\eta^{-\frac{2n^2-9n+4}{2n}C_1}.
\endaligned
\end{equation*}

Because $2n^2-9n+4>0$ for $n\geq 5$, this contradicts the
boundedness on the energy of $(\ref{eb})$. This completes the
proof.

Next, we use Proposition  $\ref{mor}$ to show that if there are
many unexceptional intervals, they must form a cascade and must
concentrate at some time $t_*$.

\begin{corollary}\label{control}
Assume that the solution $u$ is spherically symmetric. For any
interval $I \subseteq [t_-, t_+]$ and $I$ be a union of consecutive
unexceptional intervals $I_j$. Then
\begin{equation*}
\sum_{ I_j \subseteq I} \big|I_j\big|^{1/2} \lesssim
\eta^{-13C_1}\big|I\big|^{1/2},
\end{equation*}
and moreover, there exists a $j$ such that
\begin{equation*}
\big|I_j\big|\gtrsim \eta^{26C_1}\big|I\big|.
\end{equation*}
\end{corollary}

{\bf Proof: } For any unexceptional interval $I_j$,  from
H\"{o}lder inequality and Corollary $\ref{massconc2}$, we have
\begin{equation*}
\aligned \eta^{C_1}\big|I_j\big| & \lesssim \text{Mass}\big(u(t),
B(0,
\eta^{-3C_1}|I_j|^{1/2})\big) \\
& \lesssim \bigg\| \Big|
\chi\Big(\frac{x}{\eta^{-3C_1}|I_j|^{1/2}}\Big) \Big|^2 |x|^3
\bigg\|_{L^{\infty}_x} \bigg\| \Big|
\chi\Big(\frac{x}{2\eta^{-3C_1}|I_j|^{1/2}}\Big) \Big|^2
\frac{|u(t,x)|^2}{|x|^3} \bigg\|_{L^1_x} \\
& \lesssim \Big(\eta^{-3C_1}|I_j|^{1/2} \Big)^3 \int_{|x|\leq
2\eta^{-3C_1}|I_j|^{1/2}} \frac{|u(t,x)|^2}{|x|^3} dx,
\endaligned
\end{equation*}
therefore
\begin{equation*}
\aligned  \int_{|x|\leq 2\eta^{-3C_1}|I_j|^{1/2}}
\frac{|u(t,x)|^2}{|x|^3} dx \gtrsim \eta^{10C_1}\big|I_j
\big|^{-\frac12}.
\endaligned
\end{equation*}
We integrate this over each unexceptional interval $I_j$ and sum
over $j$,
\begin{equation*}
\aligned \eta^{10C_1} \sum_{ I_j \subseteq I} \big|I_j
\big|^{\frac12} & \lesssim  \sum_{ I_j \subseteq I}
\int_{I_j}\int_{|x|\leq 2\eta^{-3C_1}|I_j|^{1/2}}
\frac{|u(t,x)|^2}{|x|^3} dx\\
& \lesssim  \sum_{ I_j \subseteq I} \int_{I_j}\int_{|x|\leq
2\eta^{-3C_1}|I|^{1/2}} \frac{|u(t,x)|^2}{|x|^3} dx\\
& \lesssim   \int_{I}\int_{|x|\leq
2\eta^{-3C_1}|I|^{1/2}} \frac{|u(t,x)|^2}{|x|^3} dx\\
& \lesssim \eta^{-3C_1}|I|^{1/2}.
\endaligned
\end{equation*}
The second claim follows from the first and the fact that
\begin{equation*}
\aligned \big| I_j \big|^{1/2} \geq \big|I_j \big| \big(\sup_{I_k
\subseteq I} \big|I_k \big| \big)^{-1/2}.
\endaligned
\end{equation*}
This completes the proof.

\begin{proposition}[Interval cascade]\label{cascade}
Let $I$ be an interval tiled by finitely many intervals $I_1,
\cdots, I_N$. Suppose that for any continuous family $\big\{ I_j:
j\in \mathcal{J} \big\}$ of the unexceptional intervals, there
exists $j_* \in \mathcal{J}$ such that
\begin{equation}\label{removec}
\big| I_{j_*}\big| \geq a \big| \bigcup_{j\in \mathcal{J}} I_j \big|
\end{equation}
for some small $a>0$. Then there exist $K\geq \log(N)/\log(2a^{-1})$
distinct indices $j_1, \cdots, j_K$ such that
\begin{equation*}
\big|I_{j_1}\big| \geq 2\big|I_{j_2}\big| \geq \cdots \geq
2^{K-1}\big|I_{j_K} \big|,
\end{equation*}
and for any $t_* \in I_{j_K}$,
\begin{equation*}
\mathrm{dist}(I_{j_k}, t_*)\lesssim \frac1a \big| I_{j_k} \big|
\end{equation*}
hold for $1\leq k\leq K$.
\end{proposition}

{\bf Proof: } Here we use an algorithm in \cite{Bo98} and
\cite{Tao05} to assign a generation to each $I_j$.

\begin{figure}
\centering
\includegraphics[width=0.7\textwidth]{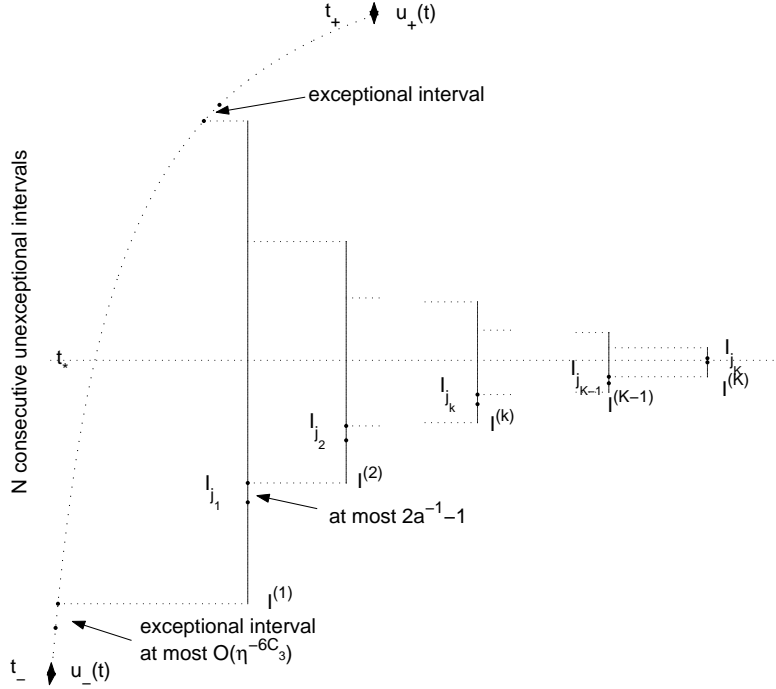}
\caption[]{ Iteration process in Proposition \ref{cascade}. }
\end{figure}
By hypothesis, $I$ contains at least one interval of length
$a|I|$. All intervals with length larger than $ a|I|/2$ belong to
the first generation. By the total measure, we see that there are
at most $2a^{-1}-1$ intervals in the first generation. Removing
there intervals from $I$ leaves at most $2a^{-1}$ gaps, which are
tiled by intervals $I_j$.

By $(\ref{removec})$ and the contradiction argument, we know that
there is not gap with length larger than $ |I|/2$.

We now apply this argument recursively to all gaps generated by the
previous iteration until every $I_j$ has been labeled with a
generation number.

Each iteration of the algorithm removes at most $2a^{-1}-1$ many
intervals and produces at most $2a^{-1}$ gaps. Suppose that there
are $N$ consecutive unexceptional intervals initially, and we
perform at most $K$ times iterations. Then the number $K$ obeys
\begin{equation*}
\aligned
 N & \leq (2a^{-1}-1) + (2a^{-1}-1)2a^{-1} + \cdots +
(2a^{-1}-1)(2a^{-1})^{K-1} \\
& \leq (2a^{-1})^K,
\endaligned
\end{equation*}
which leads to the claim $K\geq \log(N)/\log(2a^{-1})$.

Let $I^{(K)}$ be the interval obtained after $K-1$ iterations and
$I_{j_K}$ be any interval in  $I^{(K)}$. For $1\leq i \leq K-1$, let
$I^{(i)}$ be the $(i-1)$-generation gap which contains the
$I_{j_K}$, and assign the $I_{j_i}$ be any ith-generation interval
which is contained in $I^{(i)}$ (see Figure $1$). By the
construction, for any $t_* \in I_{j_K}$, we have
\begin{equation*}
\mathrm{dist}(t_*, I_{j_k}) \leq |I^{(k)}| \leq 2a^{-1}\big| I_{j_k}
\big|
\end{equation*}
for all $1 \leq k \leq K$.

\begin{proposition}[Energy non-evacuation] \label{nonevacu}
Let $I_{j_1}, \cdots, I_{j_K}$ be a disjoint family of unexceptional
intervals obeying
\begin{equation}\label{cascade2}
\big|I_{j_1}\big| \geq 2\big|I_{j_2}\big| \geq \cdots \geq
2^{K-1}\big|I_{j_K} \big|
\end{equation}
and for any $t_* \in I_{j_K}$,
\begin{equation*}
\mathrm{dist}(I_{j_k}, t_*)\lesssim \eta^{-26C_1}\big| I_{j_k} \big|
\end{equation*}
hold for $1\leq k\leq K$. Then
\begin{equation*}
K\leq \eta^{-100C_1}.
\end{equation*}
\end{proposition}

{\bf Proof: } By Corollary $\ref{massconc2}$,
\begin{equation*}
\mathrm{Mass}(u(t), B(0, \eta^{-3C_1}|I_{j_k}|^{1/2})) \gtrsim
\eta^{C_1} |I_{j_k}|
\end{equation*}
for all $t\in I_{j_k}$. By $(\ref{meR})$, we have
\begin{equation*}
\aligned \mathrm{Mass}(u(t_*), B(0, \eta^{-27C_1}|I_{j_k}|^{1/2}))
& \gtrsim\Big( \big( \eta^{C_1}
|I_{j_k}| \big)^{1/2} - \frac{\mathrm{dist}(t_*, I_{j_k})}{\eta^{-27C_1}|I_{j_k}|^{1/2} }\Big)^{2}\\
& \gtrsim \eta^{C_1} |I_{j_k}|.
\endaligned
\end{equation*}

On the other hand, from $(\ref{mer})$, we have
\begin{equation*}
\mathrm{Mass}(u(t_*), B(0, 2\eta^{C_1}|I_{j_k}|^{1/2})) \lesssim
\eta^{2C_1} |I_{j_k}|.
\end{equation*}

Define
\begin{equation*}
A(k)=\big\{ x :  \eta^{C_1}|I_{j_k}|^{1/2} \leq |x| \leq
\eta^{-27C_1}|I_{j_k}|^{1/2} \big\},
\end{equation*}
then we have
\begin{equation*}
\aligned \int_{A(k)} \big| u(t_*,x)\big|^2 dx & \gtrsim
\mathrm{Mass}(u(t), B(0, \eta^{-27C_1}|I_{j_k}|^{1/2})) -
\mathrm{Mass}(u(t), B(0, 2\eta^{C_1}|I_{j_k}|^{1/2})) \\
& \gtrsim \eta^{C_1} |I_{j_k}|.
\endaligned
\end{equation*}
By H\"{o}lder inequality, we have
\begin{equation*}
\aligned \int_{A(k)}\big| u(t_*, x) \big|^{\frac{2n}{n-2}} dx &
\gtrsim \big( \eta^{C_1}|I_{j_k}|\big)^{\frac{n}{n-2}} \big(
\eta^{-27C_1}|I_{j_k}|^{1/2} \big)^{-\frac{2n}{n-2}} \\
& \gtrsim \eta^{95C_1}
\endaligned
\end{equation*}

Choosing $M=-56C_1 \log \eta$, then we obtain by $(\ref{cascade2})$
\begin{equation*}
\aligned \eta^{-27C_1}|I_{j_{M+1}}|^{1/2} & \leq \eta^{C_1}
|I_{j_1}|^{1/2} ;\\
\eta^{-27C_1}|I_{j_{2M+1}}|^{1/2} & \leq \eta^{C_1}
|I_{j_{M+1}}|^{1/2} ;\\
& \cdots
\endaligned
\end{equation*}
Hence the annuli $A(k)$ associated to $k=1, M+1, 2M+1, \cdots, $ are
disjoint. The number of such annuli is $O(K/M)$.

Therefore from $(\ref{eb})$, we obtain
\begin{equation*}
\aligned \frac{K}{M} \eta^{95C_1} \lesssim \int_{\mathbb{R}^n}\big|
u(t_*, x) \big|^{\frac{2n}{n-2}} dx \lesssim 1.
\endaligned
\end{equation*}

That is
\begin{equation*}
K \lesssim M \eta^{-95C_1} \lesssim \eta^{-100C_1}.
\end{equation*}

We now return to the proof of Theorem $\ref{main}$. As explained at
the beginning of this section, it suffices to bound the number of
the unexceptional intervals.

Note that the number of exceptional interval is at most
$O(\eta^{-6C_3})$. We first bound the number $N$ of unexceptional
intervals that can occur consecutively.

Let us denote the union of these consecutive unexceptional
intervals by $I$. By Corollary $\ref{control}$, the hypotheses of
Proposition $\ref{cascade}$ are satisfied with $a=\eta^{26C_1}$
and so we can find a cascade of $K$ intervals and they satisfied
the hypotheses of Proposition $\ref{nonevacu}$. The bound on $K$
implies the bound on $N$, namely,
\begin{equation*}
\aligned N \lesssim(2a^{-26C_1})^{K}\thickapprox
(2\eta^{-26C_1})^{\eta^{-100C_1}}.
\endaligned
\end{equation*}

At last, since there are at most $O(\eta^{-6C_3})$ exceptional
intervals, the total number of intervals is
\begin{equation*}
J \lesssim \eta^{-6C_3} + \eta^{-6C_3}N \lesssim e^{\eta^{-200C_1}}.
\end{equation*}

This completes the proof of Theorem \ref{main}.

\textbf{Acknowledgements:} The authors were partly supported by the
NNSF of China. G. Xu wish to thank Xiaoyi Zhang for providing the
paper \cite{KiVZ} and some discussions.

\begin{center}

\end{center}
\end{document}